\documentclass[12pt,a4paper]{article}      
\usepackage{graphicx}
\usepackage{enumerate}
\usepackage{amssymb, amsmath, amsthm}

\newtheorem{theorem}{Theorem}[section]
\newtheorem{lemma}[theorem]{Lemma} 
\newtheorem{proposition}[theorem]{Proposition} 
 
\newtheorem{corollary}[theorem]{Corollary} 
 
\newtheorem{remark}[theorem]{Remark}

\begin{document}

\begin{center}
{\LARGE{\textbf{Generalized Fibonacci Numbers\\[0.4cm] and Blackwell's Renewal Theorem}}}
\end{center}




\begin{center}\vspace{.8cm}{\textbf{ S\"oren Christensen}} 

Christian-Albrechts-Universit\"at, Mathematisches Seminar, Kiel, Germany
\end{center}




\vspace{.8cm}
\textbf{Abstract:}
We investigate a connection between generalized Fibonacci numbers and renewal theory for stochastic processes. Using Blackwell's renewal theorem we find an approximation to the generalized Fibonacci numbers. With the help of error estimates in the renewal theorem we figure out an explicit representation.\\[.2cm]

\textbf{Keywords:} {generalized Fibonacci numbers,~ renewal theory}\\[.2cm]
\textbf{Mathematics Subject Classification (2000):} {11B39, 60K20}

\section{Introduction}
\label{intro}
A variety of generalizations for the Fibonacci numbers has been proposed over the last decades and were studied using different techniques. In this article we study the $d$-generalized Fibonacci numbers ($d\geq 2$) defined as follows:
\[F_n^{(d)}=0\mbox{~~~for all~}n\leq 0,~~~F_1^{(d)}=1\]
and by the recursion
\begin{equation}\label{recursion}
F_n^{(d)}=\sum_{i=1}^dF_{n-i}^{(d)}\mbox{~~~for all~}n\geq 2.
\end{equation}
Obviously the standard Fibonacci numbers come out for $d=2$. The case $d=3$ is also of special interest, the resulting sequence is called Tribonacci sequence and the first few terms for $n\geq 0$ are given by
\[0,~1,~1,~2,~4,~7,~13,~24,~44,~81,~149,...\]
The generalized Fibonacci numbers were introduced by E. P. Miles in \cite{M} and were further studied e.g. in \cite{F} and \cite{LLKS}, where Binet-type representations were given. These kind of representations were obtained using different standard techniques for difference equations such as generating functions and matrix techniques. From these techniques one expects to obtain representations of the form
\begin{equation}\label{binet}
F_n^{(d)}=\sum_{i=1}^d\lambda_iz_i^n\mbox{~~~for all~}n\geq 0,
\end{equation}
where $z_1,...,z_d$ are the (possibly complex) roots of the polynomial $x^d-x^{d-1}-...-x-1$. The coefficients $\lambda_1,...,\lambda_d$ are not so straightforward to find. \\

In this article we take a stochastic point of view to represent the generalized Fibonacci numbers via a stochastic process and use Blackwell's renewal theorem we find the asymptotic behavior of the generalized Fibonacci numbers using one special root of the above polynomial in Section \ref{sec:1}. Here it is remarkable that the right constant comes out immediately.\\
Since we wish to get an explicit representation of $F_n^{(d)}$ for all $n$ we use geometric error estimates for the renewal theorem to see that our approximation is sufficiently good in Section \ref{section_error}. In opposite to the representation in (\ref{binet}) our representations is just based on one special root of the polynomial.

\section{Approximation using renewal theory}
\label{sec:1}

First we give a combinatorial interpretation of the $d$-generalized Fibonacci numbers. $F_{n}^{(d)}$ is the numbers of possibilities to write $n-1$ as the sum of numbers in $\{1,...,d\}$:
\begin{lemma}\label{comb}
With the notation
\[\mathcal{F}_n^{(d)}=\bigcup_{m=0}^\infty\{(x_i)_{i=1}^m\in\{1,...,d\}^m:\sum_{k=1}^mx_k=n\}\mbox{~~~for all~}n\in \mathbb{Z}\]
it holds that
\[F_n^{(d)}=|\mathcal{F}_{n-1}^{(d)}|\mbox{~~~for all~}n\in\mathbb{Z}.\]
\end{lemma}
\begin{proof}
For all $n\leq 0$ it holds that $F_n^{(d)}=0=|\emptyset|=|\mathcal{F}_{n-1}^{(d)}|$ and for $n=1$ we have $F_1^{(d)}=1=|\{()\}|=|\mathcal{F}_{0}^{(d)}|$.
For $n>1$ the following recursion holds
\[|\mathcal{F}_n^{(d)}|=|\bigcup_{k=1}^d\{(x_i)_{i=1}^m\in\mathcal{F}_{n}^{(d)}:x_m=k\}|=\sum_{k=1}^m|\mathcal{F}_{n-k}^{(d)}|,\]
which proves the claim.
\end{proof}
For the use of renewal theory we embed the generalized Fibonacci numbers into a stochastic process as follows:\\
 Let $X_1,X_2,...$ be independent and identically distributed random variables with $P(X_1=i)=q^i$ for all $i=1,...,d$, where $q=q_d$ is the unique number in $(0,1)$ with 
\begin{equation}\label{q}
q+q^2+...+q^d=1
\end{equation}
and let $S_k=\sum_{i=1}^kX_i, k=1,2,...$ be the random walk generated by $X_1,X_2,...$. Furthermore let $\tau_n=\inf\{k\geq 0: S_k\geq n\}$ denote the first passage time over $n$ for $n=1,2,...$.

\begin{remark}
Formally for $d=\infty$ the random variables $X_1,X_2,...$ have a geometric distribution with parameter $1/2$. For geometric distributions renewal theory turns out to be much easier than in the general case.
\end{remark}
\begin{remark}
By multiplying equation (\ref{q}) by $1/q^d$ we see that $1/q$ is a root of the polynomial $x^d-x^{d-1}-...-x-1$.
\end{remark}

\begin{proposition}
For all $n\geq 1$ it holds that 
\[P(S_{\tau_n}=n)=q^nF_{n+1}^{(d)}.\]
\end{proposition}

\begin{proof}
For all $(x_1,...,x_m)\in\mathcal{F}^{(d)}_n$ we have
\[P(X_1=x_1,...,X_m=x_m)=\prod_{i=1}^mP(X_i=x_i)=q^{\sum_{i=1}^mx_i}=q^n,\]
hence using Lemma \ref{comb}
\[P(S_{\tau_n}=n)=\sum_{(x_i)_{i=1}^m\in\mathcal{F}^{(d)}_n}P(X_1=x_1,...,X_m=x_m)=q^n|\mathcal{F}^{(d)}_n|=q^nF_{n+1}^{(d)}.\]
\end{proof}

The proposition above states that we can represent $F_n^{(d)}$ in terms of the probability that $S_1,S_2,...$ visits the state $n$. This probability is well studied in probability theory, more specifically in renewal theory. The random variables $X_1,X_2,X_3,...$ are interpreted as renewal lifetimes and $P(S_{\tau_n}=n)$ can be seen as the probability for a renewal at time $n$. The asymptotic behavior of $P(S_{\tau_n}=n)$ for $n\rightarrow\infty$ is well understood:\\
The Blackwell renewal theorem for the lattice case (cf. e.g. \cite[Theorem I.2.2]{A}) states that 
\[P(S_{\tau_n}=n)\rightarrow\frac{1}{EX_1}=\frac{1}{\sum_{i=1}^di q^i}=\frac{(q-1)^2}{dq^{d+2}-(d+1)q^{d+1}+q}=:c_d\]
for $n\rightarrow\infty$. This leads to the following 
\begin{theorem}\label{asympt_equi}
\[F_{n}^{(d)}\approx c_dq^{-(n-1)},\]
where $\approx$ means asymptotic equivalence.
\end{theorem}

\begin{proof}
Using Lemma \ref{comb} it holds that $F_{n+1}^{(d)}q^n=P(S_{\tau_n}=n)\rightarrow c_d.$
\end{proof}
\begin{remark}
With equation (\ref{binet}) in mind it is not hard to believe that $F_{n}^{(d)}q^{n-1}$ converges to a constant, but it is not clear how to find the right constant easily. As seen above for our approach nearly no any algebra is needed.
\end{remark}

\section{Error estimation}\label{section_error}
When using the approximation $c_dq^{-(n-1)}$ for $F_n^{(d)}$ it is first of all not clear if the approximation is good or not. As an example let us consider the Tribonacci sequence
\[0,~1,~1,~2,~4,~7,~13,~24,~44,~81,~149,...\]
In this case $c_3q^{-(n-1)}$ is given by
\[0.33,~0.61,~1.13,~2.09,~3.84,~7.07,~13.01,~23.94,~44.03,~80.99,~148.98,...\]
and we see that the approximation is close to the true value; the error is always less than $1/2$, so that rounding the approximation gives the correct value.\\

In this section we estimate the error in the general case when using the approximation to the generalized Fibonacci numbers obtained in the previous section, i.e. we study
\begin{equation}\label{error}
x_n=F^{(d)}_{n}-c_dq^{-(n-1)}=q^{-(n-1)}[P(S_{\tau_{n-1}}=n-1)-c_d],~~~n\in\mathbb{Z}.
\end{equation}
Because we want to find a global bound for $|x_n|$ we especially have to show that $|P(S_{\tau_n}=n)-c_d|$ has a geometric convergence rate. The question of convergence speed in the renewal theorem has been studied intensively and we apply these results in our setting. To this end we first show that $X_1$ is "new better then used":

\begin{lemma}\label{nbu}
For all $i,j\geq 0$ with $P(X_1>i)>0$ it holds that
\[P(X_1>i+j|X_1>i)\leq P(X_1>j).\]
\end{lemma}

\begin{proof}
The proof is just straightforward calculus:
\[P(X_1>i+j)=1-P(X_1\leq i+j)=1-\sum_{k=1}^{i}q^k-\sum_{k=i+1}^{i+j}q^k\]
and
\[P(X_1>i)P(X_1>j)=(1-\sum_{k=1}^iq^k)(1-\sum_{k=1}^{j}q^k)=1-\sum_{k=1}^{i}q^k-(1-\sum_{k=1}^{i}q^k)\sum_{k=1}^jq^k,\]
therefore $P(X_1>i+j|X_1>i)\leq P(X_1>j)$ is equivalent to
\[\sum_{k=i+1}^{i+j}q^k\geq \sum_{k=1}^jq^k(1-\sum_{k=1}^{i}q^k),\]
i.e. equivalent to
\[q^i\geq 1-\sum_{k=1}^{i}q^k=\sum_{k=i+1}^{d}q^k=q^i\sum_{k=1}^{d-i}q^k\]
what is obviously true.
\end{proof}

Now we can prove the following 

\begin{theorem}\label{estimate}
\[|x_n|\leq (1-q)\left(\frac{1-q}{q}\right)^{n-1}\mbox{~~~for all~}n\geq1.\]
\end{theorem}

\begin{proof}
For the proof we use the results from \cite{B} on the rate of convergence in the renewal theorem: Using Lemma \ref{nbu} we can apply \cite[Corollary 2.1]{B} and obtain
\[|P(S_{\tau_{n-1}}=n-1)-c_d|\leq P(X_1>1)^n\mbox{~~~~for all $n\geq 1$,}\]
hence
\begin{align*}
|x_n|&=q^{-(n-1)}|P(S_{\tau_{n-1}}=n-1)-c_d|\\
&\leq q^{-(n-1)}P(X_1>1)^n=q^{-(n-1)}(1-q)^n\\
&=(1-q)\left(\frac{1-q}{q}\right)^{n-1}
\end{align*}
for each $n\geq 1$.
\end{proof}

As a result we obtain the representation of $F_n^{(d)}$ in terms of $c_d$ and $q$ as follows:

\begin{corollary}
$|x_n|<1/2$ for all $n\in\mathbb{Z}$, so that $F_n^{(d)}$ is the unique integer with minimal distance to $c_dq^{n-1}$.
\end{corollary}

\begin{proof}
 Since $q>1/2$ it holds that $(1-q)/q<1$, so by $(i)$ we have $|x_n|<1/2$ for all $n\geq 1$. Hence it remains to show $|x_n|<1/2$ for all $n\leq 0$. Since in this case $|x_n|=c_dq^{-(n-1)}$ we only have to treat $|x_0|=q/E(X_1)$. We have
\begin{align*}
\frac{1}{q}E(X_1)=\sum_{i=1}^diq^{i-1}&\geq 1+2q+q^2+...+q^{d-1}\\
&>1+q+q^2+...+q^{d-1}+q^d=2,
\end{align*}
so that $|x_0|<1/2$.
\end{proof}

\begin{remark}
We just want to mention that the previous result can also be obtained elementary without the use of renewal theory. Furthermore the representation given above is different to the usual Binet-type representation since it just uses one root of the characteristic polynomial. This has the benefit that not all (complex) roots of the polynomial have to be found.
\end{remark}

\begin{remark}
With equation (\ref{binet}) in mind we can furthermore see that $1/q$ is the only root of the equation $x^d-x^{d-1}-...-x-1$ with modulus $\geq 1$, because else $x_n=F_n^{(d)}-c_dq^{n-1}$ could not converge to 0. This is one more proof for this well known result.
\end{remark}

\end{document}